\documentclass{amsart}
\usepackage[utf8]{inputenc}
\usepackage{amsmath, amsthm, latexsym, amssymb, amscd,amsfonts}
\usepackage{tikz-cd}
\usepackage{geometry}
\usepackage{yfonts}
\usepackage{mathtools}  
\usepackage{diffcoeff}  
\usepackage[hidelinks]{hyperref}
\usepackage{amsrefs}

\renewcommand{\epsilon}{\varepsilon}
\renewcommand{\phi}{\varphi}

\newcommand{\abs}[1]{\left| #1 \right|}

\newcommand{\Gm}{\mathbb{G}_\mathrm{m}}

\renewcommand{\c}{\gamma}

\newcommand{\s}{\sigma}

\renewcommand{\P}{\mathbb{P}}
\newcommand{\Q}{\mathbb{Q}}
\newcommand{\Qbar}{\overline{\mathbb{Q}}}

\newcommand{\R}{\mathbb{R}}
\newcommand{\Z}{\mathbb{Z}}

\newcommand{\sub}{\subseteq}
\newcommand{\dif}{\setminus}

\newtheorem{thm}{Theorem}[section]

\newtheorem{prop}[thm]{Proposition}
\newtheorem{lem}[thm]{Lemma}

\title{Division of primitive Points in an abelian Variety}
\author{Francesco Ballini}
\address[F. Ballini]{Mathematical Institute, University of Oxford, Woodstock Road OX2 6GG, Oxford, United Kingdom}
\email{ballini@maths.ox.ac.uk}
\date{\today}

\begin{document}
\begin{abstract}
Let $A$ be an abelian variety defined over a number field $K$. We say that a point $P \in A(\Qbar)$ is primitive if there is no $Q \in A(\Qbar)$ defined on the field of definition of $P$ over $K$ such that $[N]Q=P$ for some positive integer $N \ge 2$. For any primitive point $P \in A(\Qbar)$, positive integer $N$ and point $Q \in A(\Qbar)$ such that $[N]Q=P$, we prove an effective lower bound on the degree of the field of definition of $Q$ over $K$ of the form $N^{\delta}$ that depends only on $A,K$ and the degree of the field of definition of $P$ over $K$. The proof is based on the estimates of the degree of torsion points by Masser.

We combine this result with a uniform version of Manin-Mumford to prove an effective Unlikely Intersections-type result: if $P \in A(\Qbar)$ is primitive, defined over a field of degree $d$ over $K$, and $X$ is a subvariety of $A$, then $X \cap [N]^{-1}P$ is contained in the weakly special part of $X$, provided $N$ is bigger than a suitable power of $d$.

As an application, we study an inverse elliptic Fermat equation, analogous to a modular Fermat equation treated by Pila.
\end{abstract}
\maketitle

\section{Introduction}

Let $A$ be an abelian variety defined over a number field $K$; we denote with $\Qbar$ an algebraic closure of $K$. We say that a point $P \in A(\Qbar)$ is \textit{primitive} if there is no $Q \in A(\Qbar)$ defined on the field of definition of $P$ over $K$ such that $[N]Q=P$ for some $N \ge 2$, where $[N]$ denotes the multiplication by $N$ on $A$. Our main result is the following:

\begin{thm}\label{thm:maindegree}
Let $A$ be an abelian variety of dimension $g$ defined over a number field $K$. There is a positive constant $c$ such that for every primitive point $P \in A(\Qbar)$, with field of definition of degree $d$ over $K$, and every point $Q \in A(\Qbar)$ such that $[N]Q=P$ for some positive integer $N \ge 2$, then the degree of the field of definition of $Q$ over $K$ is at least
\begin{equation*}
c \cdot \sqrt{d} \cdot N^{(1/4g)} / \sqrt{\log{N}}
\end{equation*}
Moreover, $c$ is effectively computable in terms of $g$, the degree of $K$ over $\Q$ and suitable equations defining $A$.
\end{thm}

Notice that the primitivity condition is necessary: one could take any $Q$ and, setting $P=[N]Q$, the bound on the degree of the field of definition of $Q$ cannot depend on $N$; this yields indeed a non-primitive $P$.

Primitivity seems to be quite common (see Section 5 for a quantitative discussion). One can construct infinitely many primitive points by choosing a number field $K'$ such that $A(K')$ contains two independent points $P_1,P_2$ (i.e. if $[N_1]P_1+[N_2]P_2$ is the origin then $N_1=N_2=0$); by the Mordell-Weil Theorem $A(K')$ is a finitely generated abelian group $G$ of rank at least 2 and hence the indivisible elements of $G$ correspond to primitive points of $A(K')$.\\

The proof is based on a combinatorial argument: if $P$ is primitive then $Q$ has a conjugate (over the field of definition of $P$), say $Q'$, such that $Q-Q'$ is a torsion point of order at least $\sqrt{N}$. One concludes using the degree estimates for torsion points, due to Masser (\cite{MASSER}).\\

Theorem \ref{thm:maindegree} is connected to Kummer Theory. If we denote by $K_N$ the field generated over $K$ by the $N$-torsion points of $A$, then the Galois group of the extension $K_N \sub K_N([N]^{-1}P)= K_N(Q)$ was studied first by Ribet in \cite{RIBET} and then by Hindry in \cite{HINDRY} (see for instance Lemme 14 and subsequent Proposition 1 for a statement similar to Theorem \ref{thm:maindegree}). While the exponent of $N$ is generally much better in these contexts, our bound is (to our knowledge) anyway not an immediate consequence of their methods, mostly because $P$ is allowed to vary in bounded extensions of $K$, not just in a fixed number field, and the constants arising from Kummer theoretic arguments have a bad dependency on $P$ (unlike the degree bounds obtained by transcendental methods, e.g. the arguments of Masser in \cite{MASSER}, which have a weaker exponent but are much more explicit). However, one can obtain Theorem \ref{thm:maindegree} in some cases, provided having an (effective) Open Image Theorem. In the case of $A$ being an elliptic curve, Davide Lombardo kindly explained to us how to deduce a better bound using the arguments of a paper by himself and Tronto \cite{LOMBARDO}.\\

We can apply the Theorem above to obtain an Unlikely Intersections-type result. We point out that the Kummer theoretic results similar to Theorem \ref{thm:maindegree} would not suffice for the bound given below.

\begin{thm}\label{thm:mainintersection}
Let $X \sub A$ be a subvariety defined over $K$ and let $X_{ws}$ be the union of the all the translates of positive dimensional abelian subvarieties of $A$ contained in $X$. Then, there is a positive constant $L$ such that, for any positive integers $d \ge 2$ and $N$ satisfying:
\begin{equation*}
N > L \cdot d^{2g} \cdot (\log{d})^{2g}
\end{equation*}
And primitive point $P \in A(\Qbar)$ defined over a field of degree $d$ over $K$, then $[N]^{-1}P$ does not intersect $X \dif X_{ws}$ (in particular, if $X$ contains no positive dimensional translate of an abelian subvariety, then $[N]^{-1}P$ does not intersect $X$). Moreover, $L$ is effectively computable in terms of $g$, the degree of $K$ over $\Q$, suitable equations defining $A$ and suitable equations defining $X$.
\end{thm}

The primitivity condition is necessary (as before, we could reseverse-engineer a suitable $P$) and the new condition on the translates of abelian subvarieties in $X$ is needed as well.

To see this, let us put ourselves in the situation $X=E \times E$, with $E$ an elliptic curve. Observe that as $X=E \times \{O\}$ (the point $O$ being the origin of $E$), then for every primitive point $p \in E(\Qbar)$ and positive integer $N$, then $(p,O)$ is primitive and $[N]^{-1}(p,O)$ intersects $X=X_{ws}$.

This is a rather special situation, since $X$ is a torsion translate. For non-torsion ones, fix a primitive non-torsion point $p \in E(\Qbar)$ and let $X=E \times \{p\}$. Observe that, for any positive integer $N$, the point $(p,[N]p)$ is primitive of bounded degree and $[N]^{-1}(p,[N]p)$ intersects $X=X_{ws}$.

The subcript of $X_{ws}$ here stands for the \textit{weakly special} part of $X$; in this case the weakly special subvarieties are the translates of abelian subvarieties (in contrast, special ones need to be torsion translates).\\

The proof of this result is based on a combination of Theorem \ref{thm:maindegree} with a uniform version of Manin-Mumford, proved independently by Hrushovski in \cite{HRUSHOVSKI} and R{\'e}mond in \cite{REMOND}. Notice that Theorem \ref{thm:mainintersection} is not an immediate consequence of Mordell-Lang-type results since $P$ is allowed to vary in all fields of a fixed degree, whose compositum is not a finitely generated field.\\

The paper is organized as follows. Section 2 is devoted to the proof of Theorem \ref{thm:maindegree} using the degree estimates of Masser. In Section 3 we prove Theorem \ref{thm:mainintersection} using uniform Manin-Mumford. Section 4 contains a treatment of an Inverse Elliptic Fermat Equation (analogous to the Modular Fermat Equation of Pila of \cite{PILAMOD}) using the methods of the previous sections. In Section 5 we briefly discuss the likeliness of a point being primitive in this context.

\section{Estimating the degree of Q}

In order to obtain a high degree for the field of definition, we will exploit torsion points with high order. Results of this kind have been proven by Masser in \cite{MASSER} and David in \cite{DAVID}. This is Theorem 4 of \cite{BERTRAND}:

\begin{thm}[Masser]\label{masserthm}
Let $A$ be an abelian variety of dimension $g$ defined over a number field $K$. Then there is a positive constant $c'$ such that the field of definition of every torsion point of order $N \ge 2$ has degree over $K$ which is at least
\begin{equation*}
c' \cdot N^{(1/g)} / (\log{N})
\end{equation*}
Moreover, $c'$ is effectively computable in terms of $g$, the degree of $K$ over $\Q$ and suitable equations defining $A$.
\end{thm}

The main ingredient for our proof is to show that if we have points $[N]Q=P$ for $P$ primitive and $N$ big, then we obtain a torsion point with big order. This is achieved using the following strategy: if $Q$ is not defined on the field of definition of $P$, then for any Galois conjugate $Q'$ we will have $[N]Q'=P$ and hence $Q-Q'$ will be a torsion point of order dividing $N$. The primitivity will ensure that such order is indeed big.\\

We sketch the situation if $N$ is prime: since $P$ is primitive, then $Q$ is not defined on the same field of $P$, hence for any conjugate $Q' \neq Q$ we have that $[N](Q-Q')=O$ (the origin). Since $N$ is prime and $Q-Q'$ is nonzero, then the degree over $K$ of $Q-Q'$, by Theorem \ref{masserthm}, is at least $c' \cdot N^{(1/g)}/\log{N}$, hence we obtain the degree of $Q$ taking the square root, as the degree of $Q-Q'$ is at most the square of the degree of $Q$.

Let us see what happens if $N$ is of the form $2M$ with $M$ prime: we have that $[2M](Q'-Q)=O$. The order can be either $2,M,2M$. In the latter two cases we are done for big $N$. If $[2](Q-Q')=O$ for every conjugate of $Q$, then the point $[2]Q$ is stable by Galois conjugation, hence it is defined over the same field of $P$, but then $P=[M][2]Q$ is not primitive.\\

In order to apply this strategy for a general $N$, we need the following combinatorial tool:

\begin{prop}\label{ordellprop}
Let $G$ be a finite abelian group and let $N$ be the minimum positive integer such that $N \cdot g=0$ for every $g \in G$. Let $S$ be a subset of $G$ with the following properties:
\begin{enumerate}
\item The lcm of the orders of the elements of $S$ is $N$;
\item For every $a,b \in S$, there exists $c \in S$ such that $a-b$ and $c$ have the same order.
\end{enumerate}
Then $S$ contains an element of order $\ge \sqrt{N}$
\end{prop}

\textit{Proof}: Let us write $G$ as $G_1 \times \ldots \times G_k$ where each $G_i$ is a group of order the power of a prime (and distinct $\abs{G_i}$ are powers of different primes). Let $N_i$ be the order of the highest order element of $G_i$. We have $N=N_1 \cdot \ldots \cdot N_k$. We denote with $p_i: G \rightarrow G_i$ the projection arising from the above product.

Let $z_1, \ldots , z_n$ be a sequence of elements of $S$ constructed as follows. $z_1$ is chosen randomly; we suppose that the order of $p_i(z_1)$ is $N_i$ precisely for $i=1, \ldots ,r_1$. Then we choose $z_2$ such that for some $i>r_1$ then the order of $p_i(z_2)$ is $N_i$ (this is ensured by condition 1); we assume that the $i>r_1$ for which the order of $p_i(z_2)$ is $N_i$ are precisely $r_1+1, \ldots , r_2$ (notice that we do not care about the orders of $p_i(z_2)$ for $i \le r_1$). In general, we choose $z_{j+1}$ so that for some $i>r_j$ the order of the projection $p_i(z_{j+1})$ is $N_i$ and we permute the $G_i$ so that the $i>r_j$ for which the order of $p_i(z_{j+1})$ is $N_i$ are exactly $r_j+1, \ldots , r_{j+1}$. We keep going until we reach the minimal $n$ for which $r_n=k$.

Summary: $z_j$ is such that $p_i(z_j)$ has order $N_i$ if $r_{j-1}+1 \le i \le r_j$, such that $p_i(z_j)$ has order strictly dividing $N_i$ if $r_j+1 \le i$ and we do not care about the order of $p_i(z_j)$ if $i \le r_{j-1}$.\\

We now consider the set $S'$ of $(x,j) \in S \times \{0, \ldots n-1\}$ such that:

\begin{center}
For $r_j+1 \le i \le k$, the product of the $N_i$ such that the order of $p_i(x)$ is $N_i$ is at least $\sqrt{N_{r_j+1} \cdot \ldots \cdot N_k}$.
\end{center}

Where we set $r_0=0$. The pair $(z_n,n-1)$ witnesses that $S$ is nonempty: in particular, for each $r_{n-1}+1 \le i \le k$ we have that $p_i(z_n)$ has order $N_i$.

We take an element of $S$ with minimal $j$. If $j=0$, then the corresponding $x$ has order at least $\sqrt{N_1 \cdot \ldots \cdots N_k}=\sqrt{N}$, so we are done.

If $j>0$, with corresponding element $x$, consider the elements $x$ and $x-z_j$. We have that if $p_i(x)$ has order $N_i$ for some $i \ge r_j+1$ then also $p_i(x-z_j)$ has order $N_i$, since the order of $p_i(z_j)$ strictly divides $N_i$.

If we have instead $r_{j-1}+1 \le i \le r_j$, we have that the order of $p_i(z_j)$ is exactly $N_i$, so that at least one of $p_i(x)$ and $p_i(x-z_j)$ has order $N_i$ (recall that the $G_i$ have a prime power order). But then, the product of the $N_i$ such that $N_i$ is the order of either $p_i(x)$ or $p_i(x-z_j)$ is precisely $N_{r_j-1}+1 \cdot \ldots \cdots N_{r_j}$; thus, either the $N_i$ which are equal to the order of $p_i(x)$ have product at least $\sqrt{N_{r_{j-1}+1} \cdot \ldots \cdots N_{r_j}}$ - but then $(x,j-1) \in S'$ - or the same property holds for $x-z_j$. That is, we have that the $N_i$ which are equal to the order of $p_i(x-z_j)$ for $r_{j-1}+1 \le i \le k$ have their product at least $\sqrt{N_{r_{j-1}+1} \cdot \ldots \cdots N_k}$. By condition 2, there is some $y \in S$ with the same order of $x-z_j$ (hence same orders of projections), so that $(y,j-1) \in S'$, contradicting the minimality of $j$. $\square$ \\

The bound $\sqrt{N}$ is probably close to optimal. One can obtain examples of $S$ consisting only of elements of order near $N^{4/7}$. We consider the configuration of the seven lines of $\P_2(\mathbb{F}_2)$ (the Fano plane), associating each point $P_i$ to a prime $p_i$. Let us take the group $G=C_{p_1} \times \dots \times C_{p_7}$, where $C_p$ denotes the cyclic group of order $p$. We choose $S$ consisting of 7 elements, one, say $a$, for each line $l$: we require $a \equiv 0 \pmod{p_i}$ if $P_i \in l$ and $a \equiv 1 \pmod {p_i}$ if $P_i \not \in l$. If the $p_i$ are big and close to each other, one obtains the bound $N^{4/7}$.\\

\textit{Proof of Theorem \ref{thm:maindegree}}: We apply Proposition \ref{ordellprop} in the following context: let $G$ be the group of $N$-torsion points of $A$. Since $[N]Q=P$, for any conjugate $Q'$ of $Q$ (with respect to a Galois automorphism which fixes $P$ and $K$), we have that $Q'-Q \in G$. Let $S$ be the set of $Q'-Q$ for the various $Q'$ obtained by Galois conjugation over the field of definition of $P$ over $K$.

First, if condition 1 was not satisfied, then there would be a proper divisor $M \mid N$ such that $[M](Q'-Q)=0$ for every conjugate $Q'$. This means that $[M]Q$ is defined on the same field as $P$ and therefore $P$ is not primitive.

Second, condition 2 is satisfied: we just need to show that for every pair of conjugates $Q_1$ and $Q_2$ of $Q$, there is a conjugate $Q_3$ such that $(Q_2-Q)-(Q_1-Q)=Q_2-Q_1$ and $Q_3-Q$ have the same order. The action of the Galois group is transitive on the various $Q'$, so take an automorphism $g$ with $g(Q_1)=Q$; just set $Q_3=g(Q_2)$ and notice that $g$ preserves orders in $G$.

From Proposition \ref{ordellprop} we conclude that for some conjugate $Q'$ we have that the order of $Q'-Q$ is at least $\sqrt{N}$. By Masser's Theorem \ref{masserthm} we obtain that $Q'-Q$ has degree (over $K$) at least $2c' \cdot N^{(1/2g)} / \log{N}$. Since $Q$ and $Q'$ have the same degree over the field of definition of $P$, we observe that:

\begin{center}
$(\text{Degree of }Q\text{ over the field of definition of }P)^2 \cdot d \ge (\text{Degree of }Q-Q'\text{ over }K)$
\end{center}
we then obtain the bound $\sqrt{2c'} \cdot \sqrt{d} \cdot N^{(1/4g)} / \sqrt{\log{N}}$ as we notice that the degree of $Q$ over $K$ is just $d$ times the degree of $Q$ over the field of definition of $P$ (the equality $[N]Q=P$ implies that the field of definition of $P$ over $K$ is a subfield of the field of definition of $Q$ over $K$). $\square$\\

The same proof has a multiplicative analogue: if $p \in \Gm(\Qbar)$ is not a perfect power in its field of definition (i.e. primitive), then every $q \in \Gm(\Qbar)$ such that $q^N=p$ has a conjugate $q'$ such that $q/q'$ is a root of unity of order at least $\sqrt{N}$. Hence, fixed an $\epsilon >0$, the degree of $q$ over $\Q$ is at least $c(\epsilon) \cdot N^{(1/4)-\epsilon}$ for some explicit $c(\epsilon)$.

\section{Uniform Manin-Mumford}

We now prove Theorem \ref{thm:mainintersection}. We use the following formulation of a Theorem of Hrushovski (Theorem 1.1.3 of \cite{HRUSHOVSKI}), which is a uniform effective Manin-Mumford-type result. See also the work of R{\'e}mond \cite{REMOND}.

\begin{thm}[Hrushovski]\label{hrushovskithm}
Let $A$ be an abelian variety defined over $K$ with a fixed embedding in a projective space. Let $Y \sub A$ be a subvariety, defined over $\Qbar$. There are a positive integer $M$ and translates $T_1, \ldots , T_M$ of abelian subvarieties of $A$ such that:
\begin{itemize}
\item Each $T_i$ is contained in $Y$
\item Each point of $Y$ which is torsion in $A$ is contained in some $T_i$
\item $M \le L' \cdot (\deg{Y})^e$
\end{itemize}
Where $L'$ and $e$ are effectively computable in terms of $g$, the degree of $K$ over $\Q$ and suitable equations defining $A$; they do not depend on $Y$.
\end{thm}

Here, the degree of $Y$ is the usual degree of algebraic geometry (i.e., if $Y$ is irreducible, the cardinality of the intersection of $Y$ and a ``generic'' linear subspace of codimension $\dim{Y}$ and, for reducible $Y$, the sum of the degrees of its irreducible components).\\

\textit{Proof of Theorem \ref{thm:mainintersection}:} We recall that $X$ here is defined over $K$. Let $P \in A(\Qbar)$ be a primitive point whose field of definition is an extension of $K$ of degree $d$. Let $Q \in A(\Qbar)$ such that $P=[N]Q$ for some positive integer $N$. Let us assume that $Q \in X \dif X_{ws}$; we recall that $X_{ws}$ is the union of all the positive dimensional translates of abelian subvarieties of $A$ contained in $X$. Notice that, if $Q'$ is a Galois conjugate of $Q$ over $K$ (and hence over the field of definition of $P$), then $Q' \in X \dif X_{ws}$. Indeed, $X_{ws}$ is stable by Galois conjugation, as if some translate of an abelian subvariety $T$ is contained in $X$, then each Galois conjugate $T'$ will be a translate of an abelian subvariety contained in $X$ as well.

We consider a translation of $X$:
\begin{center}
$X-Q=\{R \in A \mid R+Q \in X\}$
\end{center}

And notice that the degree of $X-Q$ is the same as the degree of $X$. Let $s$ be the number of conjugates of $Q$ over the field of definition of $P$; we have that $d \cdot s$ is the degree of the field of definition of $Q$ over $K$. For each of these conjugates $Q'$, we have that $Q'-Q$ is a torsion point (as $[N]Q'=P$, since the conjugation fixes $P$) lying in $X-Q$. Moreover, $Q'-Q$ cannot be contained in any translate $T$ of a positive dimensional abelian subvariety of $A$ contained in $X-Q$, for otherwise $Q' \in T+Q \sub X$, contradicting $Q' \not \in X_{ws}$. By Theorem \ref{hrushovskithm} applied on $Y=X-Q$, we have:

\begin{center}
$s \le L' \cdot (\deg{X})^e$
\end{center}
Compare with Theorem \ref{thm:maindegree}:

\begin{center}
$d \cdot s \ge c \cdot \sqrt{d} \cdot N^{(1/4g)} / \sqrt{\log{N}}$
\end{center}
Combining the two, we obtain:

\begin{center}
$N^{(1/4g)} / \sqrt{\log{N}} \le L'' \cdot \sqrt{d} \text{ } \text{ } \text{ (1)}$
\end{center}
For an explicit constant $L''$ which depends only on $X$. Suppose by contradiction that:

\begin{center}
$N > L \cdot d^{2g} \cdot (\log{d})^{2g}$
\end{center}
For some big $L$. Notice that the left hand side of (1) is increasing in $N$ as $N \ge \exp{(2g)}$, hence we have:

\begin{center}
$N^{(1/2g)} / \log{N} > L^{1/2g} \cdot d \cdot \log{d} / (\log{L}+2g\log{d}+2g\log{\log{d}}) \ge L^{1/2g} \cdot d \cdot \log{d} / (4g\log{L}\log{d}) \ge (L'')^2 \cdot d$
\end{center}
As we choose an explicit $L$, depending on $X$ only. $\square$

\section{An inverse elliptic Fermat Equation}

Pila proved in 2015 in \cite{PILAMOD} (Theorem 1.1) a modular analogue of the Fermat Equation. The proof used the Pila-Wilkie Counting Theorem (see \cite{PILAWILKIE}), coming from Model Theory, together with a strategy implemented by Pila and Zannier in \cite{PILAZANNIER}. We can prove an elliptic equivalent in the case of the \textit{inverse} Fermat equation using the arguments of the previous sections; we do not use the Pila-Wilkie Theorem, which is not effective, and our results instead are.\\

Let $E$ be an elliptic curve defined over $\Q$, given by an affine equation $y^2=x^3+ax+b$ for some $a,b \in \Q$, completed with a point at infinity obtained by homogenization, which will be the origin. We are interested in the points $P \in E(\Qbar)$ such that $x(P)$ is rational; we call these $x$-\textit{rational} points.

We say that an $x$-rational point $P$ is $x$-primitive if there is no $x$-rational $Q$ such that $P=[N]Q$ for some $N \ge 2$. Notice that $P$ being primitive implies $P$ being $x$-primitive, while the vice-versa is not true; for instance, if $a=-5$ and $b=4$, the point $P$ given by $x(P)=4,y(P)=2 \sqrt{3}$ is $x$-primitive, but it is not is not primitive since it is equal to $[2]Q$, with $x(Q)=\sqrt{3},y(Q)=1-\sqrt{3}$.\\

We can now state our result. Let $E_1$ and $E_2$ be two elliptic curves defined over $\Q$, respectively by affine equations $y_1^2=x_1^3+a_1x_1+b_1$ and $y_2^2=x_2^3+a_2x_2+b_2$, with $a_1,b_1,a_2,b_2 \in \Q$, both completed with a point at infinity obtained by homogenization, again our origin.

\begin{thm}\label{thm:invell}
There exists a positive constant $L>0$ such that, for every $N_1,N_2$ with $\max\{N_1,N_2\}>L$ and every pair of $x$-primitive $x$-rational points $P_1 \in E_1(\Qbar)$ and $P_2 \in E_2(\Qbar)$, there are no points $Q_1 \in E_1(\Qbar)$ and $Q_2 \in E_2(\Qbar)$ such that
\begin{center}
$[N_1]Q_1=P_1, [N_2]Q_2=P_1$ and $x_1(Q_1)+x_2(Q_2)=1$
\end{center}
Moreover, $L$ can be computed effectively in terms of $E_1$ and $E_2$ (namely, $a_1,b_1,a_2,b_2$).
\end{thm}

Also in this case, the $x$-primitivity is necessary: one could fix $x$-rational points $Q_1$ and $Q_2$ and construct $P_1=[N_1]Q_1$ and $P_2=[N_2]Q_2$ by reverse-engineering.

Notice that, in contrast to the Theorem 1.1 in \cite{PILAMOD}, we do not require $\max\{N_1,N_2\}$ to be prime and, moreover, our results are effective. On the other hand, such a Theorem encompasses both the direct and the inverse cases of the Fermat equation (indeed, they are the same in the modular world); the direct case here seems much more difficult, as it is in the context of $\Gm$ (see for instance this paper by Lenstra \cite{LENSTRA}).

We do want not to use general Weierstrass equations of the form $y^2=x^3+Ax^2+Bx+C$, since if $A$ is zero then the subvariety given by $x_1+x_2=1$ never contains elliptic curves. Notice that if $E_1=E_2$ the equation $x_1=x_2$ would define a union of two elliptic curves and, if $A$ is allowed to vary, one can have the same issue also for the equation $x_1+x_2=1$ (for instance, for $E_1$ given by $y_1^2=x_1^3+Ax_1^2+Bx_1+C$ and $E_2$ given by $y_2^2=-(1-x_2)^3-A(1-x_2)^2-B(1-x_2)-C$).\\

In order to exploit our previous results, we first relate $x$-primitivity with our old notion of primitivity, adapting Theorem \ref{thm:maindegree}. Given an elliptic curve $E$ defined by $y^2=x^3+ax+b$ with $a,b \in \Q$, we denote by $C_E$ the minimum positive integer such that:

\begin{center}
If $T$ is a torsion point defined either over $\Q$ or over a number field of degree 2 over $\Q$, the order of $T$ divides $C_E$.
\end{center}

Such a minimum exists, since torsion points of bounded degree are finitely many (and they are effectively computable, e.g. their height is explicitely bounded).

\begin{lem}\label{lemxprim}
Let $P \in E(\Qbar)$ be an $x$-primitive $x$-rational point and let $Q \in E(\Qbar)$ be a point such that $Q$ is defined on the field of definition of $P$. If, for some positive integer $N$, we have $[N]Q=P$, then $N$ divides $C_E$.
\end{lem}

\textit{Proof}: We use the coordinates $x$ and $y$ to denote the corresponding points (e.g. $P=(x(P),y(P))$). Suppose that $[N](x(Q),y(Q))=(x(P),y(P))$. If $P$ is rational, then $Q$ is rational as well and therefore $N = 1$ by $x$-primitivity. So suppose that $P$ is not rational, i.e. $y(P)$ is irrational.

The field of definition of $P$ has degree 2 over $\Q$, so let $\s$ be the nontrivial automorphism of the Galois group of $\Q(P)/\Q$. We have:
\begin{equation*}
[N](\s(x(Q)),\s(-y(Q)))=\s[N](x(Q),-y(Q))=\s(x(P),-y(P))=(x(P),y(P))=[N](x(Q),y(Q))
\end{equation*}

So that the point $(\s(x(Q)),\s(-y(Q)))-(x(Q),y(Q))$ is torsion of degree 1 or 2 and order $M$, with $M$ dividing $N$. But then:
\begin{equation*}
(\s(x([M]Q)),\s(-y([M]Q)))=[M](\s(x(Q)),\s(-y(Q)))=[M](x(Q),y(Q))=(x([M]Q),y([M]Q))
\end{equation*}

And this means that $[M]Q$ is $x$-rational, hence $M=N$ by $x$-primitivity of $P$. Since $M$ is the order of a torsion point of degree 1 or 2, then $M$, and hence $N$, divides $C_E$. $\square$

\begin{prop}\label{prop:maindegree}
There is a positive constant $C$ such that for every $x$-primitive $x$-rational point $P \in E(\Qbar)$ and every point $Q \in E(\Qbar)$ such that $[N]Q=P$ for some positive integer $N > C_E$, then the degree of the field of definition of $Q$ over $\Q$ is at least
\begin{equation*}
C \cdot N^{(1/4)} / \sqrt{\log{N}}
\end{equation*}
Moreover, $C$ is effectively computable in terms of $E$ (i.e. in terms of $a$ and $b$).
\end{prop}

\textit{Proof}: Let us consider all the conjugates $Q'$ of $Q$ over the field of definition of $P$. We set $G$ as the group generated by all the $Q'-Q$ for varying $Q'$; $G$ is a subgroup of the torsion points of $E$ of order $N$. As before, $S$ consists of the elements $Q'-Q$.

Let $m$ be the lcm of the orders of the elements of $S$ (and hence $G$). Notice that $m$ divides $N$, so write $N=m \cdot m'$. We have that $R=[m]Q'$ is independent of $Q'$, hence it is defined on the field of definition of $P$. We observe that $[m']R=P$, so that $m' \le C_E$ by the Lemma above. This implies that $m \ge N / C_E$.

By Proposition \ref{ordellprop} with $G$ and $S$, we obtain that there is some $Q'$ such that $Q'-Q$ is a torsion point of order at least $\sqrt{m} \ge \sqrt{N}/\sqrt{C_E}$. By the estimates of Masser (Theorem \ref{masserthm}), we have that the degree (over $\Q$) of the field of definition of $Q'-Q$ is at least $(2c'/\sqrt{C_E}) \cdot \sqrt{N}/(\log{N}-\log{C_E})$. We obtain the degree of the field of definition of $Q$ by taking the square root. $\square$\\

\textit{Proof of Theorem \ref{thm:invell}}: Let us first show that $X \sub E_1 \times E_2$, given by the equation $x_1+x_2=1$, does not contain positive dimensional translates of abelian subvarieties of $E_1 \times E_2$. We claim that $X$ is irreducible (over $\Qbar$), arguing as follows: the elliptic curve $E_1$ has function field $\Qbar(x_1,y_1)$ given by $y_1^2=x_1^3+a_1x_1+b_1$; by the equation $x_1+x_2=1$ we obtain $x_2=1-x_1$ and hence

\begin{center}
$y_2^2=(1-x_1)^3+a_2(1-x_1)+b_2$
\end{center}

We claim by contradiction that $y_2$ cannot be an element of $\Qbar(x_1,y_1)$. Notice that the right hand side of the above equality is symmetric via the involution $(x_1,y_1) \rightarrow (x_1,-y_1)$, so that $y_2$ would be either of the form $f(x_1)$ or $y_1f(x_1)$ for some rational function $f(x_1) \in \Qbar(x_1)$ (for instance, write $y_2=p(x_1)+y_1q(x_1)$). None of the above possibilities can occur, respectively because the right hand side is not a square and because $(1-x_1)^3+a_2(1-x_1)+b_2$ and $x_1^3+a_1x_1+b_1$ have at least a non-common zero (e.g. the sum of the zeroes of the first is 3, while for the second it vanishes - here the special form of our Weierstrass equations plays a role), hence with a different parity.

Let then $Y$ be a nonsingular model of the function field $\Qbar(x_1,y_1,x_2,y_2)$. The argument above implies that one can construct a dominant rational map from $Y$ to $X$ (namely, the inclusion $\Qbar(x_1,y_1) \sub \Qbar(x_1,y_1,y_2)$ gives a two-to-one morphism from $Y$ to $E_1$), hence $X$ is the closure of the image of an irreducible curve, thus irreducible.\\

Now we just need to show that $X$ is not a translate of a positive dimensional abelian suvariety of $E_1 \times E_2$; notice that $X$ passes through the origin of $E_1 \times E_2$, so, if that was the case, it would be a subgroup itself. For big $x_1$, observe that the $x$-coordinate of $[2](x_1,y_1)$ is equal to $x_1/4+O(\abs{1/x_1})$ (as long as the Weierstrass form is $y_1^2=x_1^3+a_1x_1+b_1$ - no degree two term). The same holds for $x_2$, hence for any point $P \in X$ close to the origin we have $x_1([2]P)+x_2([2]P)=1/4+O(\abs{1/x_1(P)})$, so such a curve does not give a subgroup.\\

Suppose now that we are in the setting of Theorem \ref{thm:invell} and suppose that $N_1 \ge N_2$. We consider all the conjugates of $Q_1$ and $Q_2$, call a generic pair $Q_1'$ and $Q_2'$, over the field of definition of $P_1$ and $P_2$; these are at least (referring to Proposition \ref{prop:maindegree}):

\begin{center}
$\frac{1}{4} \cdot C \cdot N_1^{(1/4)} / \sqrt{\log{N_1}}$
\end{center}

And the differences $(Q_1'-Q_1,Q_2'-Q_2) \in E_1 \times E_2$ give rise to such many torsion points on:

\begin{center}
$X-(Q_1,Q_2)=\{(R_1,R_2) \in E_1 \times E_2 \mid (R_1+Q_1,R_2+Q_2) \in X\}$.
\end{center}

As we observed, $X$ contains no positive dimensional translate of abelian subvariety of $E_1 \times E_2$, so Theorem \ref{hrushovskithm} of Hrushovski gives us an upper bound on the number of torsion points of $X-(Q_1,Q_2)$:

\begin{center}
$L' \cdot (\deg{X})^e$
\end{center}

Which is an absolute constant depending explicitely on $E_1$ and $E_2$ (namely, $a_1,b_1,a_2,b_2$). This indeed gives the desired explicit bound on $\max\{N_1,N_2\}$. $\square$

\section{Primitivity is generic}

Let us consider our elliptic curve $E: y^2=x^3+ax+b$ with $a,b \in \Q$. From the point of view of the Mordell-Weil Theorem, primitivity does not seem pervasive: the group of $\Q$-rational points $E(\Q)$ is isomorphic to a finitely generated group and indivisible points in such groups are usually a positive proportion (if the rank is at least 2), but there are cases where they are finitely many, e.g. when the group is isomorphic to $\Z$. However, such a comparison is misleading: in fact, the Weil height of the points in the Mordell-Weil groups grows steeply (as one can compare the height with a positive definite quadratic form on $\R \otimes_{\Z} E(\Q)$) and hence, while counting points of bounded height and bounded degree, it seems more likely to obtain primitive (e.g. small) elements.\\

As an example, let us count, in terms of an appropriate height, the number of the $x$-primitive $x$-rational points of $E(\Qbar)$. We refer to the Weil height, so that for instance $h(a/b)=\log{\max\{\abs{a},\abs{b}\}}$ whenever $a,b$ are relatively prime integers. We define the height of $P \in E(\Qbar)$ as $h(P)=h(x(P))$ and we recall that we can define a canonical height as $\hat{h}(P)=\lim_{n \rightarrow +\infty} h([2^n]P)/4^n$. It can be proved that there is a constant $\c$, depending only on $E$ (explicitely in terms of $a$ and $b$) such that $\abs{\hat{h}(P)-h(P)} < \c$. Such a canonical height has the nice property that $\hat{h}([N]P)=N^2 \hat{h}(P)$ for every positive integer $N$ and $P \in E(\Qbar)$.\\

Let $Z(T)$ be the number of non-torsion $x$-rational points of $E(\Qbar)$ such that their height $h$ is less than $\log{T}$. Then $Z(T)$ is roughly twice the number of rational numbers with height bounded by $T$ (we forget about the three points with $y$-coordinate equal to zero and the $x$-rational torsion points - finitely many), hence we have that there is a positive constants $\mu$ such that:
\begin{center}
$\mu^{-1} T^2 < Z(T) < \mu T^2$
\end{center}
(actually, $Z(T)/T^2$ tends to $24/\pi^2$ as $T \rightarrow +\infty$).

We define $\hat{Z}(T)$ similarly, as the number of non-torsion $x$-rational points of $E(\Qbar)$ such that their canonical height $\hat{h}$ is less than $\log{T}$. We have that:
\begin{center}
$e^{-2\c}\mu^{-1} T^2 < Z(e^{-\c}T) \le \hat{Z}(T) \le Z(e^{\c}T) < e^{2\c}\mu T^2$
\end{center}

Now notice that the number of non-torsion non-$x$-primitive points with canonical height less than $\log{T}$ is at most:

\begin{center}
$\hat{Z}(T^{1/4})+\hat{Z}(T^{1/9})+\ldots+\hat{Z}(T^{1/(\nu \log{T})}) \le T$ for big $T$
\end{center}

For some constant $\nu$ depending on a minimal nonzero canonical height of $x$-rational points (which exists). As $\hat{Z}(T) \ge e^{-\c}\mu^{-1} T^2$, then almost all of the $x$-rational points of bounded height are $x$-primitive.

One can also exclude the points of $E(\Q)$: the number of rational points with canonical height bounded by $\log{T}$ is bounded by some power $(\log{T})^{\eta}$ for big $\eta$ depending only on $E$ (this is due to the fact that the height on the Mordell-Weil group $E(\Q)$ is comparable to a positive definite quadratic form on the finite-dimesnional $\R$-vector space $E(\Q) \otimes_{\Z} \R$, i.e. points of height less than $t$ dist less than $t$ from the origin).

\section*{Acknowledgement}
We thank Jonathan Pila for suggesting us to work on the Inverse Elliptic Fermat Equation, for many useful comments, advices and countless inputs. We thank Davide Lombardo for enlightening discussions on Kummer Theory and detailed explanation of his and related results. We thank Damian R{\"o}ssler for helpful comments on a first draft of this paper and for pointing us to the Kummer theoretic viewpoint. We thank Ehud Hrushovski for comments and suggestions on a preliminary version of this work.

\end{document}